\documentclass{amsart}
\usepackage{amsthm}

\newtheorem{theorem}{Theorem}
\newtheorem{lemma}[theorem]{Lemma}
\newtheorem{corollary}[theorem]{Corollary}
\newtheorem{remark}{Remark}

\numberwithin{equation}{section}

\allowdisplaybreaks[4]

\begin{document}
\title[Approximation by subsequences of matrix transform means\dots]{Approximation by subsequences of matrix transform means of Vilenkin--Fourier series} 


\author{Istv\'an Blahota}
\address{University of Ny\'\i regyh\'aza\\ Institute of Mathematics and Computer Sciences\\ S\'ost\'oi \'ut 31/b\\ H-4400 Ny\'\i regyh\'aza\\ Hungary}

\email{blahota.istvan@nye.hu}


\author{\'Akos Gy\H ory}

\address{University of Ny\'\i regyh\'aza\\ Institute of Mathematics and Computer Sciences\\ S\'ost\'oi \'ut 31/b\\ H-4400	Ny\'\i regyh\'aza\\ Hungary}

\email{gyory.akos@nye.hu}
\address{University of Miskolc\\ Institute of Mathematics\\ Egyetem \'ut 1\\ H-3515 Miskolc\\ Hungary}

\email{akos.gyory@uni-miskolc.hu}


\begin{abstract}
In this paper we study the rate of approximation by matrix transform means of a special subsequence of partial sums of Vilenkin--Fourier series in the spaces $L_p(G_m)$ $(1\le p<\infty)$ and $C(G_m)$.
\end{abstract}






\maketitle

\section{Definitions and notations}

Let $\mathbb P$ be the set of positive integers and ${\mathbb N}:={\mathbb P}\cup \{0\}$.

Let the sequence $m:=(m_{0},m_{1},\dots )$ consist of positive integers, each at least 2. Denote by $\mathbb{Z}_{m_{n}}
:=\{0,1,\ldots,m_{n}-1\}$ the additive group of integers modulo $m_{n}$. Each factor group is equipped with the discrete topology.
The normalized Haar measure $\mu$ on ${\mathbb Z}_{m_n}$ is given
by $\mu (\{ j\})=1/m_n$, where $j\in\mathbb{Z}_{m_{n}}$. $G_{m} :=
\overset{\infty}{\underset{n=0}{{\times}}} {\mathbb Z}_{m_n}$. $G_{m}$ is called the Vilenkin group. 

If the generating sequence $m$ is bounded, then $G_m$ is called a bounded Vilenkin group; otherwise, it is called an unbounded Vilenkin group. 
In the bounded case one has $R_m:=\sup_{n\in\mathbb{N}} m_n<\infty$, while in the unbounded case $\sup_{n\in\mathbb{N}} m_n=\infty$. 

In case of $m=(2,2,\dots)$, we get $G$, the so-called Walsh group.

An element of the Vilenkin group $G_{m}$ is an infinite sequence
\[
x=(x_{0},x_{1},\dots,x_{n},\dots).
\]
Each coordinate satisfies $x_{n}\in \mathbb{Z}_{m_{n}}$ $(n\in {\mathbb N})$.

The group operation on $G_{m}$ is the coordinate-wise addition (denoted by $+$),
the normalized Haar measure $\mu $ is the product measure and the topology is the product topology.
$m$-adic intervals are defined in the usual way 
\[
I_{0}(x):=G_{m},\ I_{n}(x):=\{y\in G_{m}:y=(x_{0},\dots,x_{n-1},y_{n},y_{n+1},\dots)\}
\]
for $x\in G_{m}, n\in {\mathbb P}$. They form a base for the neighbourhoods of $G_{m}$.
Let $0:=(0:i\in
{\mathbb N})\in G_{m}$ denote the null element of $G_{m}$.

We define the so-called generalized number system based on $m$
in the following way: 
\[
M_{0}:=1,\ M_{n+1}:=m_{n}M_{n}\quad  (n\in \mathbb{N}).
\]
Then every $n\in \mathbb{N}$ has a unique representation of the form $n=\sum_{k=0}^{\infty }n_{k}M_{k},$ where $n_{k}\in \mathbb{Z}_{m_{k}}\ (k\in \mathbb{N})$ 
and only a finite number of $n_{k}$ differ from zero. For a given $n\in {\mathbb P}$, the order of $n$ is defined by $|n|:=\max\{ j\in {\mathbb N}: n_j\neq 0\}$. So, it is  a non-negative integer number such that $M_{|n|}\leq n <M_{|n|+1}$.

Let $L_{p}(G_{m})$ denote the usual Lebesgue spaces on $G_{m}$ (with the
corresponding norm $\Vert .\Vert_p$), where $1\leq p<\infty$. 
For $f\in C(G_m)$ we use the supremum norm $\|f\|_\infty:=\sup_{x\in G_m}|f(x)|$.

Next, we define the modulus of continuity in $L_{p}(G_{m}), 1\leq p\leq\infty,$ of a function $f\in L_{p}(G_{m})$ by
\[
\omega_p(f,\delta):=\sup_{\lvert t\rvert<\delta}\Vert f(\cdot+t)-f(\cdot)\Vert_p,\quad \delta>0,
\]
with the notation
\[
\lvert x\rvert:=\sum_{k=0}^\infty \frac{x_k}{M_{k+1}} \quad \textrm{for all }x\in G_{m}.
\]
 When $p=\infty$ we interpret $\omega_\infty$ and $\|\cdot\|_\infty$ on $C(G_m)$.
 
The Lipschitz classes in $L_p(G_{m})$ (for each $\alpha>0$) are defined as
\[
\textrm{Lip}(\alpha,p,G_{m}):=\{ f\in L_{p}(G_{m}): \omega_p(f,\delta)=O(\delta^\alpha )
\textrm{ as } \delta\to 0\} .
\]

We introduce some concepts of the Vilenkin--Fourier analysis. 
The generalized Rademacher functions are
defined as
\[
r_{k}(x):=\exp\left(\frac{2\pi\imath x_{k}}{m_k}\right) \ (x\in G_{m},\ k\in {\mathbb N}).
\]
Vilenkin functions for $n\in\mathbb{N}$ are
\[
\psi_{n}(x):=\prod_{k=0}^{\infty }r_{k}^{n_{k}}(x).
\]

Let ${\mathcal P}_n$ be the collection of Vilenkin polynomials of order less than $n$, that is, functions of the form
\[
P(x)=\sum_{k=0}^{n-1}a_k\psi_k(x),
\]
where $n\in\mathbb{P}$ and $(a_k)$ is a sequence of complex numbers.

The $n$th Fourier-coefficient, the $n$th partial sum of the Fourier series and the $n$th Dirichlet kernel are defined by
\[
\hat f(n):=
\int_{G_{m}} f\bar \psi_{n} d\mu,
\
S_{n}(f):=\sum_{k=0}^{n-1}\hat f(k)\psi_{k},
\ D_{n}:=
\sum_{k=0}^{n-1}\psi_{k},\ D_0:=0.
\]

Fej\'er kernels are defined as the arithmetic means of Dirichlet kernels, that is, 
\[
K_{n}:=\frac{1}{n}\sum_{k=1}^{n}D_{k}.
\] 

Let $T:=(t_{k,n})_{1\le k\le n<\infty}$ be an infinite upper triangular matrix. Let us define the $n$th matrix transform mean (or linear mean) determined by the matrix $T$ 
\[
\sigma_{n}^T(f):=\sum_{k=1}^{n}t_{k,n}S_{k}(f),
\]
where $(t_{k,n}: 1\leq k\leq n, \ k\in\mathbb{P})$ be a finite sequence of non-negative numbers for each $n\in\mathbb{P}$. 

It is a common generalization of the partial sum of Fourier series and several well-known (e.g. Fej\'er, Ces\`aro, N\"orlund, weighted, Riesz) means.

The $n$th matrix transform kernel ($T$ kernel) is defined by
\[
K_{n}^{T}:=\sum_{k=1}^{n}t_{k,n}D_{k}.
\]
It is easily seen that
\[
\sigma_{n}^T(f;x)=\int_{G_{m}}f(u)K_{n}^{T}(x-u)d\mu(u).
\]
This equality (and its analogous versions for special means) shows the necessity of observing kernel functions.

We introduce the notation
$\Delta t_{k,n}:=t_{k,n}-t_{k+1,n},$ where $1\leq k\leq n$ and $t_{n+1,n}:=0$.

\section{Historical overview}

Let $(q_k:k\geq 0)$ be a sequence of non-negative numbers. The $n$th N\"orlund mean of the Vilenkin--Fourier series is  defined by
\[
\frac{1}{Q_n}\sum_{k=1}^{n}q_{n-k}S_{k}(f;x),
\]
where $Q_n:=\sum_{k=0}^{n-1}q_k$ $(n\geq 1)$. It is  always assumed that $q_0>0$ and 
\[
\lim_{n\to\infty}Q_n=\infty.
\]
In this case, the N\"orlund method generated by $(q_k)$ is regular (see \cite{MS,Zyg}) if and only if 
\[
\lim_{n\to\infty}\frac{q_{n-1}}{Q_n}=0.
\]

In the paper \cite{MS}, the rate of the approximation by N\"orlund means of Walsh--Fourier series of a function $f$ in $L^p(G)$ and in $C(G)$ (in particular, in $\textrm{Lip}(\alpha ,p,G),$ where $\alpha>0$ and $1\leq p\leq\infty$) was studied.
As special cases M\'oricz and Siddiqi obtained the earlier results given by
Yano \cite{Y2}, Jastrebova \cite{J} and Skvortsov \cite{Sk2} on the rate of the approximation by Ces\`aro means. 
In 2008, Fridli, Manchanda and Siddiqi generalized the result of M\'oricz and Siddiqi for homogeneous Banach spaces and dyadic Hardy spaces \cite{FSM}. More recently Baramidze, Goginava, K. Nagy, Memi\'c, Persson, Tephnadze, Wall and the first author presented  some results with respect to this topic \cite{T1,GN2,T2,T3}. See \cite{TE,W1}, as well.

Let $(p_k:k\geq 1)$ be a sequence of nonnegative numbers.   The $n$th weighted mean of Vilenkin--Fourier series is defined by
\[
\frac{1}{P_n}\sum_{k=1}^{n}p_{k}S_{k} (f;x),
\]
where $P_n:=\sum_{k=1}^{n}p_k$ $(n\geq 1).$ In particular case we get Vilenkin--Fej\'er means (for all $k$ set $p_k=1$).
It is always assumed that $p_1>0$ and 
\[
\lim_{n\to\infty}P_n=\infty,
\]
which is the condition for regularity.

In the paper \cite{MR}, the rate of the approximation by weighted means of Walsh--Fourier series of a function in $L^p(G)$ and in $C(G)$ (in particular, in $\textrm{Lip}(\alpha ,p,G),$ where $\alpha>0$ and $1\leq p\leq\infty$) was studied. As special cases M\'oricz and Rhoades obtained the earlier results given by 
Yano \cite{Y2}, Jastrebova \cite{J} on the rate of the approximation by Walsh--Ces\`aro means. A common generalization of these two results of M\'oricz and Siddiqi \cite{MS} and M\'oricz and Rhoades \cite{MR} was given by the authors in the paper \cite{BN2}.

Matrix transform means are common generalizations of several well-known summation methods. It follows by simple consideration that the N\"orlund means, the weighted means, as well as the Fej\'er (or the $(C,1)$) and the $(C,\alpha)$ means are special cases of the matrix transform summation method introduced above. For matrix transform means with respect to trigonometric system see e.g. the results of Chandra \cite{Cha} and Leindler \cite{Le2}. For the  Walsh system, see Blyumin \cite{Bly}. 

This paper is motivated by the work of M\'oricz, Siddiqi \cite{MS} on Walsh--N\"orlund mean method and by the results of M\'oricz, Rhoades \cite{MR} on Walsh weighted mean method. It is important to note that some ideas are coming from the paper of Chripk\'o \cite{C1}. The main aim is to investigate the rate of the approximation of matrix transform in terms of modulus of continuity  under some general conditions. The main theorem (Theorem \ref{thm-main}) gives a kind of common generalization of the two results of M\'oricz, Siddiqi on N\"orlund means \cite{MS} and M\'oricz, Rhoades on weighted means \cite{MR} to special subsequences. Moreover, we extend the results from the Walsh system to the (bounded) Vilenkin system. See also \cite{B}.

Avdispahi\'c and Pepi\'c proved some results also for Vilenkin system in the paper \cite{ref1}.

Finally, we remark that a closely related theorem for the Vilenkin system was recently proved by Blahota and K. Nagy in \cite{BN5}. 
Our estimates concern the special subsequence of indices $M_n$, in contrast to the general-index setting of their theorem.
For this subsequence our approach permits weaker assumptions on the matrix coefficients; in addition, we explicitly compute the dependence of the constant $c$ on the parameters $R_{m}$ and $c_{m}$, which was left implicit in the general setting of their work.

\section{Auxiliary results}	

The $M_n$th Dirichlet kernel has a closed form 
\begin{equation}\label{dir}
	D_{M_n}(x)=\begin{cases}
		0, & \textrm{if }x\not\in I_{n}(0),\\
		M_n, & \textrm{if } x\in I_{n}(0).
	\end{cases}
\end{equation}

It is well-known (see e.g. \cite{AVDZ}) that if $G_m$ is a bounded Vilenkin group, then the $L^1(G_m)$-norms of the Fej\'er kernels are uniformly bounded. 
More precisely, there exists a constant $c_m>0$, depending only on $m$, such that
\[
\|K_n\|_{1}\le c_m \qquad (n\in\mathbb P),
\]
moreover,
\[
c_m:=\sup_{n\in\mathbb P}\|K_n\|_{1}<\infty.
\]

\begin{lemma}\cite{BN5}\label{MS-Vilenkin} Let $n\in\mathbb{N},\ q\in\mathbb{P}\ ,1\leq q\leq m_n-1$ and $g\in {\mathcal P}_{M_n}$, $f\in L^p(G_{m})$ ($1\leq p< \infty$) or $f\in C(G_{m})$. Then 
	\[
	\left\Vert\int_{G_{m}} r_n^q(u) g(u) (f(\cdot-u)-f(\cdot))d\mu(u) \right\Vert_p\leq m_n \Vert g\Vert_1\omega_p\left(f,\frac{1}{M_n}\right)
	\]
	holds (in the case $f\in C(G_m)$ we replace $p$ by $\infty$).
\end{lemma}

\begin{lemma}\label{decomp} Let $n\in\mathbb{N}$, then we have 
	\begin{align*}
		K_{M_n}^{T}=&\sum_{j=0}^{n-1}\sum_{l=1}^{m_{j}-1}\sum_{s=1}^{l-1}\sum_{k=0}^{M_j-1}t_{lM_j +k,M_{n}}r_{j}^{s}D_{M_j}
		+\sum_{j=0}^{n-1}	\sum_{q=M_j}^{M_{j+1}-1}t_{q,M_{n}}D_{M_j}\\
		&+\sum_{j=1}^{n-1}\sum_{l=1}^{m_{j}-1}r_{j}^l\sum_{k=1}^{M_j-2}\Delta t_{lM_j+k,M_{n}}kK_{k}\\
		&+\sum_{j=1}^{n-1}\sum_{l=1}^{m_{j}-1}r_{j}^l
		t_{(l+1)M_j-1,M_{n}}(M_j-1)K_{M_j-1}\\
		&+t_{M_{n},M_{n}}D_{M_{n}}
		=:\sum_{i=1}^{5} K_{i,n}.
	\end{align*}
\end{lemma}
\begin{proof}
	Let us set $0\leq k<M_{j}$ and $0<l <m_j$, then 
	\begin{equation}\label{first}
		D_{lM_j+k}=\sum_{s=0}^{l-1}\sum_{i=0}^{M_{j}-1}\psi_{sM_j+i}+\sum_{i=0}^{k-1}\psi_{lM_j+i}
		=\sum_{s=0}^{l-1}r_{j}^{s}D_{M_j}+r_{j}^lD_{k}.
	\end{equation}	
	
	We write 
	\[
	K_{M_{n}}^{T}=\sum_{j=0}^{n-1}\sum_{v=M_{j}}^{M_{j+1}-1}t_{v,M_{n}}D_{v}
	+t_{M_{n},M_{n}}D_{M_{n}}
	=:K^{1}_{n}+K^{2}_{n}.
	\]
	Now, we apply equality (\ref{first}) for the expression $K_n^{1}$. We get 
	\begin{align*}
		K_{n}^{1}=&\sum_{j=0}^{n-1}\sum_{l=1}^{m_{j}-1}\sum_{k=0}^{M_j-1}t_{lM_j +k,M_{n}}D_{lM_j+k}\\
		=&\sum_{j=0}^{n-1}\sum_{l=1}^{m_{j}-1}\sum_{k=0}^{M_j-1}t_{lM_j +k,M_{n}}\sum_{s=0}^{l-1}r_{j}^{s}D_{M_j}
		+\sum_{j=0}^{n-1}\sum_{l=1}^{m_{j}-1}\sum_{k=0}^{M_j-1}t_{lM_j +k,M_{n}}r_{j}^lD_{k}\\
		=&	\sum_{j=0}^{n-1}\sum_{l=1}^{m_{j}-1}\sum_{s=1}^{l-1}\sum_{k=0}^{M_j-1}t_{lM_j +k,M_{n}}r_{j}^{s}D_{M_j}
		+\sum_{j=0}^{n-1}	\sum_{q=M_j}^{M_{j+1}-1}t_{q,M_{n}}D_{M_j}\\
		&+\sum_{j=0}^{n-1}\sum_{l=1}^{m_{j}-1}r_{j}^l\sum_{k=0}^{M_j-1}t_{lM_j +k,M_{n}}D_{k}
		=:K^{1,1}_{n}+K^{1,2}_{n}+K^{1,3}_{n}.
	\end{align*}
	
	Now, we use Abel's transform for the expression $K^{1,3}_{n}$. We have 
	\begin{align*}
		K_{n}^{1,3}
		=&\sum_{j=0}^{n-1}\sum_{l=1}^{m_{j}-1}r_{j}^l\sum_{k=1}^{M_j-1}t_{lM_j+k,M_{n}}D_{k}\\
		=&\sum_{j=1}^{n-1}\sum_{l=1}^{m_{j}-1}r_{j}^l\left(\sum_{k=1}^{M_j-2}\Delta t_{lM_j+k,M_{n}}kK_{k }+
		t_{(l+1)M_j-1,M_{n}}(M_j-1)K_{M_j-1}\right).
	\end{align*}	
	
	It completes the proof of Lemma \ref{decomp}.
\end{proof}

\section{Main Theorem}

\begin{theorem}\label{thm-main}
	Let $f\in C(G_{m})$ or $f\in L^p(G_{m})\ (1 \leq p < \infty)$. For every $n\in\mathbb{P}$, let $(t_{k,M_{n}}: 1\leq k\leq M_{n})$ be a finite sequence of non-negative numbers such that 
	\begin{equation}\label{c1}
		\sum_{k=1}^{M_{n}}t_{k,M_{n}}=1
	\end{equation}
	is satisfied. Assume that $G_m$ is a bounded Vilenkin group.
	
	a) If the finite sequence  $(t_{k,M_{n}}: 1\leq k \leq M_{n})$ is non-decreasing for a fixed $n$, then
	\[
	\Vert\sigma^{T}_{M_{n}}(f)-f\Vert_{p}\leq c\sum_{j=0}^{n-1}M_j\sum_{l=1}^{m_j-1}t_{(l+1)M_j-1,M_{n}} \omega_{p}\left(f,\frac{1}{M_{j}}\right)+\omega_{p}\left(f,\frac{1}{M_{n}}\right)
	\]
	holds (for $f\in C(G_{m})$ we change $p$ by $\infty$).
	
	b) If the finite sequence $(t_{k,M_{n}}: 1\leq k \leq M_{n})$  is  non-increasing for a  fixed $n$, then 
	\[
	\Vert\sigma^{T}_{M_{n}}(f)-f\Vert_{p}\leq c\sum_{j=0}^{n-1}M_{j} 
	\sum_{l=1}^{m_j-1} t_{lM_j,M_{n}}
	\omega_{p}\left(f,\frac{1}{M_{j}}\right)
	+\frac{1}{M_{n}}\omega_{p}\left(f,\frac{1}{M_{n}}\right)
	\]
	holds (for $f\in C(G_{m})$ we change $p$ by $\infty$).
	
	Here $c=R_{m}^{2}+2c_{m}R_{m}+1$ explicitly.
\end{theorem}

\begin{proof} For $C(G_{m})$ the proof is similar, even simpler.
	Let us set $f\in L^p(G_{m})$. Keeping in mind Equation \eqref{c1}, we use 
	Lemma \ref{decomp} and the usual Minkowski inequality
	
	\begin{align*}
		\Vert\sigma^{T}_{M_{n}}(f)-f\Vert_{p}=&\left(\int_{G_{m}}|\sigma^{T}_{M_{n}}(f;x)-f(x)|^{p}
		d\mu(x)\right)^{\frac{1}{p}}\\
		=&\left(\int_{G_{m}}\left|\int_{G_{m}}K_{M_{n}}^{T}(u)(f(x-u)-f(x))d\mu(u)\right|^{p}
		d\mu(x)\right)^{\frac{1}{p}}\\
		\leq&\sum_{i=1}^{5}\left(\int_{G_{m}}\left|\int_{G_{m}}K_{i,n}(u)(f(x-u)-f(x))d\mu(u)\right|^{p}
		d\mu(x)\right)^{\frac{1}{p}}	\\	
		=&\sum_{i=1}^{5}I_{i,n}.
	\end{align*}
	
	For the expression $I_{1,n}$ we apply Lemma \ref{MS-Vilenkin} and Equality \eqref{dir}. We have
	\begin{align*}
		I_{1,n}\leq& \sum_{j=0}^{n-1}\sum_{l=1}^{m_{j}-1}\sum_{s=1}^{l-1}\sum_{k=0}^{M_j-1}t_{lM_j +k,M_{n}}\left(\int_{G_{m}}\left|\int_{G_{m}}r_{j}^{s}(u)\right.\right.\\
		&\times\left.\left.D_{M_j}(u)
		(f(x-u)-f(x))d\mu(u)\right|^{p}
		d\mu(x)\right)^{\frac{1}{p}}\\	
		\leq&\sum_{j=0}^{n-1}m_j\sum_{l=1}^{m_{j}-1}(l-1)\sum_{k=0}^{M_j-1}t_{lM_j +k,M_{n}}
		\omega_{p}\left(f,1/M_{j}\right)\\
		\leq&\sum_{j=0}^{n-1}m_j^2\sum_{l=1}^{m_{j}-1}\sum_{k=0}^{M_j-1}t_{lM_j +k,M_{n}}
		\omega_{p}\left(f,1/M_{j}\right).
	\end{align*}
	In case a), we write
	\[
	I_{1,n}	\leq R_{m}^{2}\sum_{j=0}^{n-1} M_j  \sum_{l=1}^{m_{j}-1} t_{(l+1)M_j-1,M_{n}}\omega_{p}\left(f,1/M_{j}\right).
	\]
	In case b), we get 
	\[
	I_{1,n}\leq R_{m}^{2}\sum_{j=0}^{n-1}M_j\sum_{l=1}^{m_{j}-1} t_{lM_j,M_{n}}	\omega_{p}\left(f,1/M_{j}\right).
	\]
	
	Using generalized Minkowski inequality and Equality \eqref{dir}
	for the expression $I_{2,n}$, 
	we obtain
	\begin{align*}
		I_{2,n}&\leq \sum_{j=0}^{n-1}\sum_{q=M_{j}}^{M_{j+1}-1} t_{q,M_{n}} 
		\int_{G_{m}}D_{M_{j}}(u)\left(\int_{G_{m}}\left|f(x-u)-f(x)\right|^{p}d\mu(x)\right)^{\frac{1}{p}}d\mu(u)\nonumber\\
		&\leq \sum_{j=0}^{n-1}\sum_{k=0}^{M_{j+1}-M_{j}-1} t_{M_{j}+k,M_{n}} 
		\omega_{p}\left(f,1/M_{j}\right)\\
		&=\sum_{j=0}^{n-1}\sum_{l=1}^{m_j-1}\sum_{r=0}^{M_j-1} t_{lM_j+r,M_n}
		\omega_{p}\left(f,1/M_{j}\right).
	\end{align*}
	
	We  get in case a)
	\begin{align*}
		I_{2,n}&\leq\sum_{j=0}^{n-1}\sum_{l=1}^{m_j-1} \sum_{r=0}^{M_j-1} t_{(l+1)M_j-1,M_n}\omega_{p}\left(f,1/M_{j}\right)\\
		&=\sum_{j=0}^{n-1} M_j\sum_{l=1}^{m_j-1}t_{(l+1)M_j-1,M_{n}}\omega_{p}\left(f,1/M_{j}\right).
	\end{align*}
	
	We similarly calculate in case b) 
	\begin{align*}
		I_{2,n}&\leq\sum_{j=0}^{n-1}
		\sum_{l=1}^{m_j-1}\sum_{r=0}^{M_j-1} t_{lM_j,M_n}		
		\omega_{p}\left(f,1/M_{j}\right)\\
		&=\sum_{j=0}^{n-1}
		M_j\sum_{l=1}^{m_j-1}t_{lM_j,M_{n}} 
		\omega_{p}\left(f,1/M_{j}\right).
	\end{align*}
	
	Now, applying Lemma \ref{MS-Vilenkin} we get 
	\begin{align*}
		I_{3,n}
		\leq&\sum_{j=1}^{n-1}\sum_{l=1}^{m_{j}-1}\sum_{k=1}^{M_j-2}\left|\Delta t_{lM_j+k,M_{n}}\right|k\\
		&\times\left(\int_{G_{m}}\left|\int_{G_{m}}r_{j}^l(u)K_{k}(u)(f(x-u)-f(x))d\mu(u)\right|^{p}d\mu(x)\right)^{\frac{1}{p}}\\
		\leq&
		\sum_{j=1}^{n-1}m_j\sum_{l=1}^{m_j-1}\sum_{k=1}^{M_{j}-2}\left|\Delta t_{lM_{j}+k,M_{n}}\right|k \omega_{p}\left(f,1/M_j\right)\|K_{k}\|_{1}\\
		\leq &
		c_{m}R_{m}\sum_{j=1}^{n-1}\sum_{l=1}^{m_j-1}\sum_{k=1}^{M_{j}-2}\left|\Delta t_{lM_{j}+k,M_{n}}\right|k \omega_{p}\left(f,1/M_{j}\right).
	\end{align*}
	
	We write in case a) 
	\begin{align*}
		\sum_{k=1}^{M_{j}-2}\left|\Delta t_{lM_{j}+k,M_{n}}\right|k
		=&\sum_{k=1}^{M_{j}-2}
		(t_{lM_j+k+1,M_{n}}-t_{lM_j+k,M_{n}} )k\\
		=& (M_j-2)t_{(l+1)M_{j}-1,M_{n}}-\sum_{k=1}^{M_{j}-2}t_{lM_j+k,M_{n}}\\
		\leq&   
		M_j t_{(l+1)M_{j}-1,M_{n}}
	\end{align*}
	and 
	\[
	I_{3,n}\leq c_{m}R_{m}\sum_{j=1}^{n-1}M_j\sum_{l=1}^{m_j-1} t_{(l+1)M_j-1,M_{n}}\omega_{p}\left(f,1/M_{j}\right).
	\]	
	We have in case b) 
	\begin{align*}
		\sum_{k=1}^{M_{j}-2}\left|\Delta t_{lM_{j}+k,M_{n}}\right|k
		&= \sum_{k=1}^{M_{j}-2}t_{lM_j+k,M_{n}}-(M_j-2)t_{(l+1)M_{j}-1,M_{n}}\\
		&\leq\sum_{k=1}^{M_{j}-2}t_{lM_j+k,M_{n}}
		\leq M_j t_{lM_j,M_{n}}
	\end{align*}
	and 
	\[
	I_{3,n}
	\leq c_{m}R_{m}\sum_{j=1}^{n-1}M_j\sum_{l=1}^{m_j-1} t_{lM_j,M_{n}}\omega_{p}\left(f,1/M_{j}\right).
	\]
	
	Now, we estimate the expression $I_{4,n}$. Lemma \ref{MS-Vilenkin} yields
	\begin{align*}
		I_{4,n}\leq&
		\sum_{j=1}^{n-1}M_j\sum_{l=1}^{m_{j}-1}t_{(l+1)M_j-1,M_{n}}\\ &\times\left(\int_{G_{m}}\left|\int_{G_{m}}r_{j}^l(u)K_{M_j-1}(u)(f(x-u)-f(x))d\mu(u)\right|^{p}d\mu(x)\right)^{\frac{1}{p}}\\
		\leq&
		\sum_{j=1}^{n-1}m_{j}M_{j}\sum_{l=1}^{m_{j}-1}t_{(l+1)M_j-1,M_{n}} 
		\omega_{p}\left(f,1/M_{j}\right)\|K_{M_{j}-1}\|_{1}\nonumber\\
		\leq&
		c_{m}R_{m}\sum_{j=1}^{n-1}M_{j}\sum_{l=1}^{m_{j}-1}t_{(l+1)M_j-1,M_{n}} 
		\omega_{p}\left(f,1/M_{j}\right).
	\end{align*}
	So, we are ready in case a).

	We have in case b) 
	\[
	I_{4,n}\leq c_{m}R_{m}\sum_{j=1}^{n-1}M_{j}\sum_{l=1}^{m_j-1}t_{lM_{j},M_{n}}
	\omega_{p}\left(f,1/M_{j}\right).
	\]

	For estimating $I_{5,n}$ let us use generalized Minkowski inequality and Equality \eqref{dir} again.
	\begin{align*}
		I_{5,n}&\leq t_{M_{n},M_{n}} 
		\int_{G_{m}}D_{M_{n}}(u)\left(\int_{G_{m}}\left|f(x-u)-f(x)\right|^{p}d\mu(x)\right)^{\frac{1}{p}}d\mu(u)\nonumber\\
		&\leq t_{M_{n},M_{n}} 
		\omega_{p}\left(f,1/M_{n}\right).
	\end{align*}
	In case a) we use $t_{M_{n},M_{n}}\leq 1$, while in case b) $t_{M_{n},M_{n}}\leq 1/M_{n}$.

	This completes the proof of our Theorem \ref{thm-main}.
\end{proof}

\begin{remark}
	We mention that assuming \eqref{c1} is natural, because many well-known (see e.g. Fej\'er, Ces\`aro, N\"orlund, weighted, Riesz) means satisfy it and this equality is a part of regularity conditions \cite[p.~74]{Zyg}.
\end{remark}
\begin{remark}
	Comparing this result with similar theorems in \cite{BN2,BN5,BN6,BNT} where we considered the expression $\Vert\sigma_{n}^T(f)-f\Vert_{p}$, we found that for this special mean in case a) we do not need to assume the analogue of the condition $t_{n,n}=O\left(n^{-1}\right)$, that is $t_{M_{n},M_{n}}=O\left(1/M_{n}\right)$.
\end{remark}
\begin{corollary}\label{s}
	Let
	\[
	t_{k,M_{n}}:=\begin{cases} 
		0, &\textrm{if } 1\leq k<M_{n}, \\
		1, &\textrm{if } k=M_{n}.
	\end{cases}
	\]
	Then for every $f\in L_p(G_m)$, $1\le p<\infty$, and for $f\in C(G_m)$ if $p=\infty$ we get
	$\sigma_{M_{n}}^T(f)=S_{M_{n}}(f)$, therefore
	\[
	\Vert S_{M_{n}}(f)-f\Vert_{p}\leq\omega_{p}\left(f,1/M_{n}\right),
	\] 	
	and if $f\in \textrm{Lip}(\alpha,p,G_{m})$, then 
	\[
	\Vert S_{M_{n}}(f)-f\Vert_{p}=O(M_{n}^{-\alpha}).
	\] 
\end{corollary}

\end{document}